\documentclass[12pt]{amsart}
\usepackage{amsmath}
\usepackage{amstext}
\usepackage{amsfonts}
\usepackage{amsthm}
\usepackage{amscd}
\numberwithin{equation}{section}

\swapnumbers
\theoremstyle{plain}
\newtheorem{thm}{Theorem}[section]
\newtheorem{lem}[thm]{Lemma}
\newtheorem{prop}[thm]{Proposition}

\theoremstyle{remark}

\def\cO{{\mathcal O}}

\def\map#1{\ \smash{\mathop{\longrightarrow}\limits^{#1}}\ }

\def\ZZ{\mathbb{Z}}

\def\FF{\mathbb{F}}

\def\deg{\mathrm{deg} \:}

\def\ker{\mathrm{ker} \:}

\def\rk{\mathrm{rk} \:}

\def\End{\mathrm{End}}
\def\Hom{\mathrm{Hom}}

\def\lra{\longrightarrow}

\def\ra{\rightarrow}

\def\lms{\longmapsto}

\def\MM{\mathcal{M}}

\begin{document}
\title{Semistability of Frobenius direct images over curves} 
\author{Vikram B. Mehta} 
\address{Tata Institute of Fundamental Research \\ Homi Bhabha Road \\ Mumbai 400 005 \\ India}
\email{vikram@math.tifr.res.in}

\author{Christian Pauly}
\address{D\'epartement de Math\'ematiques \\ Universit\'e de Montpellier II - Case Courrier 051 \\ Place Eug\`ene Bataillon \\ 34095 Montpellier Cedex 5 \\ France}
\email{pauly@math.univ-montp2.fr}

\begin{abstract}
Let $X$ be a smooth projective curve of genus $g \geq 2$ defined over an
algebraically closed field $k$ of characteristic $p>0$. Given a semistable vector
bundle $E$ over $X$, we show that its direct image $F_*E$ under the Frobenius map $F$
of $X$ is again semistable. We deduce a numerical characterization of the
stable rank-$p$ vector bundles $F_*L$, where $L$ is a line bundle
over $X$.
\end{abstract}

\maketitle 

\section{Introduction}

Let $X$ be a smooth projective curve of genus $g \geq 2$ defined over an
algebraically closed 
field $k$ of characteristic $p>0$ and let 
$F : X \ra X_1$ be the relative $k$-linear Frobenius map. 
It is by now a well-established fact that
on any curve $X$ there exist semistable vector bundles $E$ such that their pull-back
under the Frobenius map $F^*E$ is not semistable \cite{LP}, \cite{LasP}. In order to control the
degree of instability of the bundle $F^*E$, one is naturally lead 
(through adjunction $\Hom_{\cO_{X}}(F^* E, E') =
\Hom_{\cO_{X_1}}(E, F_* E'))$ to ask whether semistability is preserved by direct image under the
Frobenius map. The answer is (somewhat surprisingly) yes. In this note we show the following result.

\begin{thm}
Assume that $g \geq 2$. If $E$ is a semistable
vector bundle over $X$ (of any degree), then $F_*E$ is also semistable.
\end{thm}

Unfortunately we do not know whether also stability is preserved by  direct image under Frobenius.
It has been shown that $F_*L$ is stable for a line bundle $L$ (\cite{LP} Proposition 1.2) 
and that in small characteristics the bundle $F_*E$ is stable for any stable bundle
$E$ of small rank \cite{JRXY}. 
The main ingredient of the proof is Faltings' cohomological criterion of semistability. We also
need the fact that the generalized Verschiebung $V$, defined as the rational map from the moduli
space $\MM_{X_1}(r)$ of semistable rank-$r$ vector bundles over $X_1$ with fixed trivial determinant to the moduli space $\MM_X(r)$ induced
by pull-back under the relative Frobenius map $F$,
$$V_r: \MM_{X_1}(r) \dashrightarrow \MM_X(r), \qquad E \lms F^* E $$
is dominant for large $r$. We actually show a stronger statement for large $r$.

\begin{prop}
If $l \geq g(p-1) +1$ and $l$ prime, then the generalized Verschiebung $V_l$ is
generically \'etale for any curve $X$. In particular $V_l$ is separable and dominant.
\end{prop}

As an application of Theorem 1.1 we obtain an upper bound of the rational invariant $\nu$
of a vector bundle $E$, defined as 
$$ \nu(E) := \mu_{max}(F^*E) - \mu_{min}(F^*E), $$
where $\mu_{max}$ (resp. $\mu_{min}$) denotes the slope of the first (resp. last) piece
in the Harder-Narasimhan filtration of $F^* E$. 

\begin{prop}
For any semistable rank-$r$ vector bundle $E$ 
$$\nu(E) \leq  min((r-1)(2g-2), (p-1)(2g-2)).$$
\end{prop} 

We note that the  inequality $\nu(E) \leq (r-1)(2g-2)$ was proved in \cite{SB} Corollary 2 
and  in  \cite{S} 
Theorem 3.1. We suspect that the relationship 
between both inequalities comes from the conjectural  fact that 
the length (=number of pieces) of the
Harder-Narasimhan filtration of $F^*E$ is at most $p$ for semistable $E$.

\bigskip

Finally we show that direct images of line bundles under Frobenius are 
characterized by maximality of the invariant $\nu$.

\begin{prop}
Let $E$ be a stable rank-$p$ vector bundle over $X$. Then the following statements are
equivalent.

\begin{enumerate}
\item There exists a line bundle $L$ such that $E = F_*L$.

\item $\nu(E) = (p-1)(2g-2).$
\end{enumerate}
\end{prop}
We do not know whether the analogue of this proposition remains true for higher rank.

\section{Reduction to the case $\mu(E) = g-1$.}

In this section we show that it is enough to prove Theorem 1.1 for semistable
vector bundles $E$ with slope $\mu(E) = g-1$.

\bigskip

Let $E$ be a semistable vector bundle over $X$ of rank $r$ and let $s$ be the integer 
defined by the equality
$$ \mu(E) = g-1 + \frac{s}{r}.$$
Applying the Grothendieck-Riemann-Roch theorem to the Frobenius map $F: X \ra X_1$, 
we obtain
$$ \mu(F_*E) = g-1 + \frac{s}{pr}.$$
Let $\pi: \tilde{X} \ra X$ be a connected \'etale covering of degree $n$ and let 
$\pi_1 : \tilde{X}_1 \ra X_1$ denote its twist by the Frobenius of $k$ (see \cite{R}
section 4). The diagram
\begin{eqnarray} \label{diag1}
\begin{CD}
\tilde{X} @>F>> \tilde{X}_1 \\
@V\pi VV   @VV \pi_1 V \\
X @>F>> X_1
\end{CD}
\end{eqnarray}
is Cartesian and we have an isomorphism
$$\pi_1^*(F_* E) \cong F_*(\pi^* E).$$
Since semistability is preserved under pull-back by a 
separable morphism of curves, we see that $\pi^* E$ is semistable.
Moreover if $F_*(\pi^* E)$ is semistable, then $F_* E$ is also semistable.

\bigskip

Let $L$ be a degree $d$ line bundle over $\tilde{X}_1$. The projection formula 
$$ F_*(\pi^* E \otimes F^*L) = F_*(\pi^* E) \otimes L$$
shows that semistability of $F_*(\pi^* E)$ is
equivalent to semistability of $F_*(\pi^* E \otimes F^*L)$.

\bigskip
Let $\tilde{g}$ denote the genus of $\tilde{X}$. By the Riemann-Hurwitz formula
$ \tilde{g} - 1 = n(g-1)$. We compute
$$ \mu(\pi^* E \otimes F^*L) = n(g-1) + n\frac{s}{r} + pd = \tilde{g} -1 +n\frac{s}{r}
+pd,$$
which gives
$$ \mu(F_*(\pi^* E \otimes F^*L)) = \tilde{g} -1 +n\frac{s}{pr} +d.$$

\begin{lem}
For any integer $m$ there exists a connected \'etale covering
$ \pi: \tilde{X} \ra X$
of degree $n = p^k m$ for some $k$.
\end{lem}

\begin{proof}
If the $p$-rank of $X$ is nonzero, the statement is clear. If the
$p$-rank is zero, we know by Corollaire 4.3.4 \cite{R} that there
exist connected \'etale coverings $Y \ra X$ of degree $p^t$ for infinitely many
integers $t$ (more precisely for all $t$ of the form $(l-1)(g-1)$ where $l$ is a large
prime). Now we decompose $m=p^s u$ with $p$ not dividing $u$. We then take a covering $Y
\ra X$ of degree $p^t$ with $t \geq s$ and a covering $\tilde{X} \ra Y$ of degree
$u$. 
\end{proof}

Now the lemma applied to the integer $m = pr$ shows existence of a connected \'etale
covering $\pi : \tilde{X} \ra X$ of degree $n= p^km$. Hence $n \frac{s}{pr}$ is an integer and we
can take $d$ such that  $ n\frac{s}{pr} +d = 0$.

\bigskip

To summarize, we have shown that for any semistable $E$ over $X$ there exists a
covering $\pi : \tilde{X} \ra X$ and a line bundle $L$ over $\tilde{X}_1$ such that
the vector bundle $\tilde{E} := \pi^* E \otimes F^* L$ is semistable with 
$\mu(\tilde{E}) = \tilde{g} - 1$
and such that semistability of $F_*\tilde{E}$ implies semistability of $F_*E$.

\section{Proof of Theorem 1.1}
 
In order to prove semistability of $F_*E$ we shall use the cohomological 
criterion of semistability due to Faltings \cite{F}. 

\begin{prop} [\cite{L} Th\'eor\`eme 2.4 ]
Let $E$ be a rank-$r$ vector bundle over $X$ with $\mu(E) = g-1$ and $l$ an
integer $> \frac{r^2}{4}(g-1)$.  Then $E$ is semistable 
if and only if there exists a rank-$l$ vector bundle $G$ with trivial
determinant such that
$$ h^0(X,E \otimes G)= h^1 (X, E \otimes G) = 0.$$ 
\end{prop}

Moreover if the previous condition holds for one bundle $G$, it holds for
a general bundle by upper semicontinuity of the function $G \mapsto h^0(X, E \otimes
G)$.

\bigskip

{\bf Remark.} The proof of this proposition (see \cite{L} section 2.4) works over any
algebraically closed field $k$.

\bigskip

By Proposition 1.2 (proved in section 4) we know that $V_l$ is dominant 
when $l$ is a large prime number. Hence 
a general vector bundle $G \in \MM_X(l)$ is of the form $F^*G'$ for some
$G' \in \MM_{X_1}(l)$. Consider a semistable $E$ with $\mu(E) = g-1$. Then by 
Proposition 3.1 $h^0(X, E \otimes G) =0$ for general $G \in \MM_X(l)$. Assuming
$G$ general, we can write $G = F^* G'$ and we obtain
by adjunction 
$$ h^0(X,E \otimes F^* G') = h^0(X_1, F_*E \otimes G') = 0.$$
This shows that $F_*E$ is semistable by Proposition 3.1.

\section{Proof of Proposition 1.2}

According to \cite{MS} section 2 it will be enough to
prove the existence of a stable vector bundle $E \in \MM_{X_1}(l)$
satisfying $F^* E$ stable and
$$ h^0(X_1, B \otimes \End_0(E)) = 0,$$
because the vector space $H^0(X_1, B \otimes \End_0(E))$ can be identified with the
kernel of the differential of $V_l$ at the point $E \in \MM_{X_1}(l)$. Here 
$B$ denotes the sheaf of locally exact differentials over $X_1$ (see \cite{R}
section 4).

\bigskip

Let $l \not= p$ be a prime number and let $\alpha \in JX_1[l] \cong JX[l]$ be a
nonzero $l$-torsion point. We denote by  
$$ \pi : \tilde{X} \ra X \qquad \text{and} \qquad \pi_1 : \tilde{X}_1 \ra X_1 $$
the associated cyclic \'etale cover of $X$ and $X_1$ and by $\sigma$ a generator
of the Galois group $\mathrm{Gal}(\tilde{X}/X) = \ZZ/l\ZZ$. 
We recall that the kernel of the Norm map
$$ \mathrm{Nm} : J\tilde{X} \lra JX$$
has $l$ connected components and we denote by
$$i:  P := \mathrm{Prym}(\tilde{X}/X) \subset J\tilde{X} $$
the associated Prym variety, i.e., the connected component containing the origin. 
Then we have an isogeny
$$ \pi^* \times i:  JX \times P \lra J\tilde{X} $$
and taking direct image under $\pi$ induces a morphism
$$ P \lra \MM_X(l),\qquad  L \lms \pi_*L.$$
Similarly we define the Prym variety $P_1 \subset JX_1$ and the morphism
$P_1 \ra \MM_{X_1}(l)$ (obtained by twisting with the Frobenius of $k$).
Note that $\pi_{1*}L$ is semistable for any $L \in P_1$ and stable for 
general $L \in P_1$ (see e.g. \cite{B}). Since $F^*(\pi_{1*}L) \cong
\pi_*(F^* L)$ --- see diagram \eqref{diag1} --- and since $F^*$ induces the Verschiebung 
$V_P : P_1 \ra P$, which is surjective, 
we obtain that $\pi_{1*}L$ and $F^*(\pi_{1*}L)$
are stable for general $L \in P_1$.  

\bigskip

Therefore Proposition 1.2 will immediately follow from the next Proposition.

\begin{prop}
If $l \geq g(p-1)+1$ then there exists a cyclic degree $l$ \'etale cover
$\pi_1: \tilde{X}_1 \ra X_1$ with the property that
$$ h^0(X_1, B \otimes \End_0(\pi_{1*}L)) = 0$$
for general $L \in P_1$.
\end{prop}

\begin{proof}
By relative duality for the \'etale map $\pi_1$ we have $(\pi_{1*} L)^* \cong 
\pi_{1*} L^{-1}$. Therefore
$$ \End(\pi_{1*}L) \cong \pi_{1*} L \otimes \pi_{1*} L^{-1} \cong 
\pi_{1*} \left( L^{-1} \otimes \pi_1^* \pi_{1*} L  \right)$$
by the projection formula. Moreover since $\pi_1$ is Galois \'etale we have a
direct sum decomposition
$$ \pi_1^* \pi_{1*} L \cong \oplus_{i=0}^{l-1} (\sigma^i)^* L.$$ 
Putting these isomorphisms together we find that
\begin{eqnarray*}
H^0(X_1, B \otimes \End(\pi_{1*}L)) & = & H^0(X_1, B \otimes \pi_{1*}\left(\oplus_{i=0}^{l-1} 
L^{-1} \otimes (\sigma^i)^* L  \right) \\ 
  & = & \oplus_{i=0}^{l-1} H^0(X_1, B \otimes \pi_{1*} ( L^{-1} \otimes (\sigma^i)^* L)) \\
  & = & H^0(X_1,B \otimes \pi_{1*} \cO_{\tilde{X}_1}) \oplus
   \oplus_{i=1}^{l-1} H^0(X_1, B \otimes \pi_{1*} ( L^{-1} \otimes (\sigma^i)^* L)).
\end{eqnarray*}

Moreover $\pi_* \cO_{\tilde{X}_1} = \oplus_{i=0}^{l-1} \alpha^i$, which implies that
\begin{equation} \label{hb}
H^0(X_1, B \otimes \End_0(\pi_{1*}L)) = \oplus_{i=1}^{l-1} H^0(X_1,B \otimes \alpha^i)
\oplus \oplus_{i=1}^{l-1} H^0(X_1, B \otimes \pi_{1*} ( L^{-1} \otimes (\sigma^i)^* L)).
\end{equation}
Let us denote for $i=1,\dots,l-1$  by $\phi_i$ the isogeny 
$$ \phi_i : P_1 \lra P_1, \qquad L \lms L^{-1} \otimes (\sigma^i)^* L.$$
Since the function  $L \mapsto h^0(X_1, B \otimes \End_0(\pi_{1*}L))$ is upper semicontinuous, it will be enough to
show the existence of a cover $\pi_1: \tilde{X}_1 \ra X_1$ satisfying
\begin{enumerate}
\item for $i= 1,\dots,l-1,  \ \ \ h^0(X_1,B \otimes \alpha^i) =0$ (or equivalently, $P$ is an ordinary abelian variety).
\item for $M$ general in $P$, \ \  $h^0(X_1, B \otimes \pi_{1*}M) =0$.
\end{enumerate}
Note that these two conditions implie that the vector space \eqref{hb} equals $\{ 0 \}$ for general
$L \in P_1$, because the  $\phi_i$'s are surjective.

\bigskip

We recall that $\ker(\pi_1^*: JX_1 \ra J\tilde{X}_1) =  \langle \alpha \rangle \cong \ZZ/l\ZZ$ and that
$$ P_1[l] = P_1 \cap \pi_1^* (JX_1) \cong \alpha^\perp/ \langle \alpha \rangle $$
where $\alpha^\perp = \{ \beta \in JX_1[l] \ \ \text{with} \ \ \omega(\alpha, \beta) =1 \}$ and
$\omega: JX_1[l] \times JX_1[l] \ra \mu_l$ denotes the symplectic Weil form.
Consider a $\beta \in \alpha^\perp \setminus \langle \alpha \rangle$. Then 
$\pi_1^* \beta \in P_1[l]$ and
$$\pi_{1*} \pi_1^* \beta = \oplus_{i=0}^{l=1} \beta \otimes  \alpha^i.$$ 
Again by upper semicontinuity of the function $M \mapsto h^0(X_1,B \otimes \pi_{1*}M)$ one observes that the conditions (1) and (2)
are satisfied because of the following lemma (take $M = \pi_1^* \beta$).

\begin{lem}
If $l \geq g(p-1)+1$ then there exists a pair $(\alpha, \beta) \in JX_1[l] \times JX_1[l]$ satisfying
\begin{enumerate}
\item $\alpha \not = 0$ and $\beta \in \alpha^\perp \setminus  \langle \alpha \rangle$,
\item for $i=1,\dots,l-1 \ \ h^0(X_1,B \otimes \alpha^i) =0$,
\item for $i=0,\dots,l-1 \ \ h^0(X_1,B \otimes \beta \otimes \alpha^i) =0$.
\end{enumerate}
\end{lem}

\begin{proof}
We adapt the proof of \cite{R} Lemme 4.3.5. We denote by $\FF_l$ the
finite field $\ZZ/l\ZZ$. Then there exists a symplectic isomorphism
$JX_1[l] \cong \FF^g_l \times \FF^g_l$, where the latter space is
endowed with the standard symplectic form. Note that composition is
written multiplicatively in $JX_1[l]$ and additively in $\FF_l^{2g}$. 
A quick computation shows
that the number of isotropic $2$-planes in $\FF^g_l \times \FF^g_l$ equals
$$ N(l) = \frac{(l^{2g}-1)(l^{2g-2}-1)}{(l^2-1)(l-1)}.$$
Let $\Theta_B \subset JX_1$ denote the theta divisor associated to $B$. Then by \cite{R}
Lemma 4.3.5 the cardinality $A(l)$ of the finite set 
$\Sigma(l):= JX_1[l] \cap \Theta_B$ satisfies
$$ A(l) \leq l^{2g-2}g(p-1).$$
Suppose that there exists an isotropic $2$-plane $\Pi \subset \FF^g_l \times
\FF^g_l$ which contains $\leq l-2$ points of $\Sigma(l)$. Then we can 
find a pair $(\alpha,\beta)$ satisfying the 3 properties of the Lemma as
follows: any nonzero point $x \in \Pi$ determines a line (=$\FF_l$-vector
space of dimension $1$). Since a line contains $l-1$ nonzero points, we
obtain at most $(l-1)(l-2)$ nonzero points lying on lines generated by
$\Sigma(l) \cap \Pi$. Since $(l-1)(l-2) < l^2-1$ there exists
a nonzero $\alpha$ in the complement of these lines. Now we note that there
are $l-1$ affine lines parallel to the line generated by $\alpha$ and the
$l$ points on any of these affine lines are of the form $\beta \alpha^i$ for
$i=0,\dots,l-1$ for some $\beta \in \alpha^\perp \setminus \langle \alpha \rangle$.
The points $\Sigma(l) \cap \Pi$ lie on at most $l-2$ such affine lines, hence
there exists at least one affine line parallel to $\langle \alpha \rangle$ avoiding $\Sigma(l)$.
This gives $\beta$.

\bigskip

Finally let us suppose that any isotropic $2$-plane contains $\geq l-1$ points
of $\Sigma(l)$. Then we will arrive at a contradiction as follows: we
introduce the set 
$$ S= \{ (x,\Pi) \ |  \ x \in \Pi \cap \Sigma(l) \ \ \text{and} \ \ 
\Pi  \ \text{isotropic 2-plane} \}.$$
with cardinality $|S|$. Then by our assumption we have 
\begin{equation} \label{in1}
|S| \geq (l-1)N(l).
\end{equation}
On the other hand, since any nonzero $x \in \FF_l^g \times \FF_l^g$ is 
contained in $\frac{l^{2g-2}-1}{l-1}$ isotropic $2$-planes, we obtain
\begin{equation} \label{in2}
|S| \leq \frac{l^{2g-2}-1}{l-1} A(l).
\end{equation}
Putting \eqref{in1} and \eqref{in2} together, we obtain
$$ A(l) \geq \frac{l^{2g} -1}{l+1}.$$
But this contradicts the inequality $A(l) \leq l^{2g-2} g(p-1)$ if 
$l \geq g(p-1) +1$. 
\end{proof}

This completes the proof of Proposition 4.1.
\end{proof}

\bigskip

{\bf Remark.}  It has been shown \cite{O} Theorem A.6 that $V_r$ is dominant for any
rank $r$ and any curve $X$, by using a versal deformation of a direct sum a $r$ line bundles.

\bigskip

{\bf Remark.} We note that $V_r$ is not separable when $p$ divides the rank $r$ and $X$ is non-ordinary.
In that case the Zariski tangent space at a stable bundle $E \in \MM_{X_1}(r)$ identifies with
the quotient $H^1(X_1, \End_0(E))/ \langle e \rangle$ where $e$ denotes the nonzero extension class
of $\End_0(E)$ by $\cO_{X_1}$ given by $\End(E)$. Then the inclusion of homotheties $\cO_{X_1} \hookrightarrow \End_0(E)$ induces 
an inclusion $H^1(X_1,\cO_{X_1}) \subset H^1(X_1, \End_0(E))/ \langle e \rangle$ and the restriction of the differential of $V_r$ at the point $E$ to $H^1(X_1,\cO_{X_1})$ coincides with the non-injective Hasse-Witt map.

\section{Proof of Proposition 1.3}

Since we already know that $\nu(E) \leq (r-1)(2g-2)$ (\cite{SB}, \cite{S}) it suffices
to show that $\nu(E) \leq (p-1)(2g-2)$.

We consider the quotient $F^* E \ra Q$ with minimal slope, i.e., $\mu(Q) = \mu_{min}(F^*E)$
and $Q$ semistable.
By adjunction we obtain a nonzero morphism $E \ra F_*(Q)$, from which we deduce
(using Theorem 1.1) that 
$$\mu(E) \leq \mu(F_*Q) = \frac{1}{p} \left(\mu_{min}(F^*E) +(p-1)(g-1) \right)$$
hence
$$\mu(F^*E) \leq \mu_{min}(F^*E) +(p-1)(g-1).$$
Similarly we consider the subbundle $S \hookrightarrow F^*E$ with maximal slope,i.e.,
$\mu(S) = \mu_{max}(F^*E)$ and $S$ semistable. Taking the dual and proceeding as above,
we obtain that
$$\mu(F^*E) \geq \mu_{max}(F^*E)  - (p-1)(g-1).$$
Now we combine both inequalities and we are done.

\bigskip

{\bf Remark.} We note that the inequality of Proposition 1.3 is sharp. The maximum $(p-1)(2g-2)$ is
obtained for the bundles $E = F_* E'$ (see \cite{JRXY} Theorem 5.3).

\section{Characterization of direct images}

Consider a line bundle $L$ over $X$. Then the direct image $F_*L$ is stable 
(\cite{LP} Proposition 1.2) and the Harder-Narasimhan filtration of $F^* F_*L$ 
is of the form (see \cite{JRXY})
$$ 0 = V_0  \subset V_1 \subset \ldots \subset V_{p-1} \subset V_p = F^* F_*L,\qquad
\text{with} \qquad V_i/V_{i-1} \cong L \otimes \omega_X^{p-i}.$$
In particular $\nu(F_*L) = (p-1)(2g-2)$. In this section 
we will show a converse statement.

\bigskip
\noindent
More generally let $E$ be a stable rank-$rp$ vector bundle with 
$\mu(E) = g-1 + \frac{d}{rp}$ for some
integer $d$ and satisfying
\begin{enumerate}
\item the Harder-Narasimhan filtration of $F^*E$ has $l$ terms.
\item $\nu(E) = (p-1)(2g-2)$.
\end{enumerate}

\bigskip

{\bf Questions.} Do we have $l\leq p$? Is $E$ of the form $E = F_*G$ for some rank-$r$ vector bundle $G$? We will give a positive answer in the case $r=1$ (Proposition 6.1).

\bigskip
Let us denote the Harder-Narasimhan filtration by
$$ 0 = V_0  \subset V_1 \subset \ldots \subset V_{l-1} \subset V_l = F^* E,\qquad
V_i/V_{i-1} = M_i.$$
satisfying the inequalities
$$ \mu_{max}(F^*E) = \mu(M_1) > \mu(M_2) > \ldots > \mu(M_l) = \mu_{min}(F^*E).$$
The quotient $F^* E \ra M_l$ gives via adjunction a nonzero map $E \ra F_*M_l$.
Since $F_*M_l$ is semistable, we obtain that $\mu(E) \leq \mu(F_*M_l)$. This implies
that $\mu(M_l) \geq g-1 + \frac{d}{r}$. Similarly taking the dual of the inclusion $M_1 \subset
F^*E$ gives a map
$F^* (E^*) \ra M_1^*$ and by adjunction $E^* \ra F_*(M_1^*)$. Let us denote 
$\mu(M_1^*) = g-1 + \delta$, so that $\mu(F_*(M_1^*)) = g-1 + \frac{\delta}{p}$.
Because of semistability of $F_*(M_1^*)$,
we obtain $-(g- 1 + \frac{d}{rp}) = \mu(E^*) \leq  \mu(F^*(M_1^*))$, hence
$\delta  \geq -2p(g-1) - \frac{d}{r}$. This implies that $\mu(M_1) \leq  (2p-1)(g-1)
+ \frac{d}{r}$. Combining
this inequality with $\mu(M_l) \geq g-1 + \frac{d}{r}$ and the assumption $\mu(M_1) -\mu(M_l) 
= (p-1)(2g-2)$, we obtain that
$$ \mu(M_1) = (2p-1)(g-1) + \frac{d}{r}, \qquad \mu(M_l) = g-1 + \frac{d}{r}.$$

Let us denote by $r_i$ the rank of the semistable bundle $M_i$. We have the equality
\begin{equation} \label{sri}
\sum_{i=1}^l r_i = rp.
\end{equation}
Since $E$ is stable and $F_*(M_l)$ is semistable and since these bundles have
the same slope, we deduce that $r_l \geq r$. Similarly we obtain that
$r_1 \geq r$.

\bigskip
\noindent
Note that it is enough to show that $r_l  = r$. Since 
$E$ is stable and $F_*M_l$ semistable and since the two 
bundles have the same slope and rank, they will be isomorphic.

\bigskip
\noindent
We introduce the integers for $i= 1,\ldots, l-1$
$$ \delta_i =  \mu(M_{i+1}) - \mu(M_i) + 2(g-1) = \mu(M_{i+1} \otimes \omega) -
\mu(M_i).$$
Then we have the equality
\begin{equation} \label{sdel2}
 \sum_{i=1}^{l-1} \delta_i = \mu(M_l) - \mu(M_1) +2(l-1)(g-1)  = 2(l-p)(g-1).
\end{equation}
We note that  if
$\delta_i <0$, then $\Hom(M_i,M_{i+1} \otimes \omega) = 0$.

\begin{prop}
Let $E$ be stable rank-$p$ vector bundle with $\mu(E) = g-1 + \frac{d}{p}$ and
$\nu(E) = (p-1)(2g-2)$. Then
$E = F_*L$ for some line bundle $L$ of degree $g-1+d$.
\end{prop}

\begin{proof}
Let us first show that $l=p$. We suppose that $l<p$. 
Then $\sum_{i=1}^{l-1} \delta_i = 2(l-p)(g-1) < 0$ so that there exists a $k \leq
l-1$ such that $\delta_k < 0$.
We may choose
$k$ minimal, i.e.,$\delta_i \geq 0$ for $i <k$. Then we have
\begin{equation}
\mu(M_k) >  \mu(M_i) + 2(g-1)  \qquad \text{for} \ i > k.
\end{equation}
We recall  that $\mu(M_i) \leq \mu(M_{k+1})$ for $i >k$. The Harder-Narasimhan filtration of 
$V_k$ is given by the first $k$ terms of the Harder-Narasimhan filtration of $F^*E$. Hence
$\mu_{min}(V_k) = \mu(M_k)$. 

Consider now the canonical connection $\nabla$ on $F^*E$ and its first
fundamental form
$$ \phi_k : V_k \hookrightarrow F^* E \map{\nabla} F^* E 
\otimes \omega_X \lra (F^*E/V_k)\otimes \omega_X.$$
Since $\mu_{min}(V_k) > \mu (M_i \otimes \omega)$ for
$i>k$ we obtain $\phi_k = 0$. Hence $\nabla$ preserves $V_k$ and since $\nabla$ has
zero $p$-curvature, there exists a subbundle $E_k \subset E$ such that $F^*E_k = V_k$.

We now evaluate $\mu(E_k)$. By assumption $\delta_i \geq 0$ for $i <k$. Hence 
$$\mu(M_i) \geq \mu(M_1) - 2(i-1)(g-1) \qquad \text{for} \   i \leq k,$$
which implies that
$$ \deg(V_k)  =  \sum_{i=1}^{k} r_i \mu(M_i) \geq \rk (V_k) \mu(M_1) -
2(g-1) \sum_{i=1}^k r_i(i-1).$$
Hence we obtain
$$ p \mu(E_k) = \mu(V_k) \geq \mu(M_1) -2(g-1) C,$$
where $C$ denotes the fraction $\frac{\sum_{i=1}^k r_i(i-1)}{\rk(V_k)}$. We
will prove in a moment that $C \leq \frac{p-1}{2}$, so that we obtain by
substitution
$$ p \mu(E_k) \geq (2p-1)(g-1) + d  - (g-1)(p-1) = p (g-1) + d = p \mu(E),$$
contradicting stability of $E$. Now let us show that $C \leq \frac{p-1}{2}$ or
equivalently 
$$ \sum_{i=1}^k ir_i \leq \frac{p+1}{2} \sum_{i=1}^k r_i.$$
But that is obvious if $k \leq \frac{p-1}{2}$. Now if $k > \frac{p-1}{2}$ we note that
passing from $E$ to $E^*$ reverses the order of the $\delta_i$'s, so that the index
$k^*$ for $E^*$ satisfies $k^* \leq \frac{p-1}{2}$. This proves that $l=p$.

\bigskip
\noindent
Because of \eqref{sri} we obtain $r_i = 1$ for all $i$ and therefore $E = F_*M_p$. 

\end{proof}

\section{Stability of $F_*E$?}

Is stability also preserved by $F_*$?

\bigskip

We show the following result in that direction.

\begin{prop}
Let $E$ be a stable vector bundle over $X$. Then $F_*E$ is simple.
\end{prop}

\begin{proof}
Using relative duality  $(F_*E)^* \cong F_*(E^* \otimes \omega_X^{1-p})$ we obtain
$$ H^0(X_1, \End(F_*E)) = H^0(X, F^*F_* E \otimes E^* \otimes \omega_X^{1-p}).$$
Moreover the Harder-Narasimhan filtration of $F^* F_*E$ 
is of the form (see \cite{JRXY})
$$ 0 = V_0  \subset V_1 \subset \ldots \subset V_{p-1} \subset V_p = F^* F_*E,\qquad
\text{with} \qquad V_i/V_{i-1} \cong E \otimes \omega_X^{p-i}.$$
We deduce that
$$H^0(X,F^*F_*E \otimes E^* \otimes \omega_X^{1-p}) = H^0(X,V_1 \otimes E^* \otimes
\omega_X^{1-p}) = H^0(X, \End(E)),$$
and we are done.
\end{proof}

\end{document}